\documentstyle[amscd,amssymb,verbatim,12pt]{amsart}

\theoremstyle{plain}
\newtheorem{Thm}[subsection]{Theorem}
\newtheorem{Cor}[subsection]{Corollary}
\newtheorem{Lem}[subsection]{Lemma}
\newtheorem{Prop}[subsection]{Proposition}

\theoremstyle{definition}
\newtheorem{Def}[subsection]{Definition}

\theoremstyle{remark}

\newtheorem{Rem}[subsection]{Remark}

\numberwithin{equation}{section}
\renewcommand{\rm}{\normalshape}

\newif\ifShowLabels
\ShowLabelstrue
\newdimen\theight
\def\TeXref#1{%
	\leavevmode\vadjust{\setbox0=\hbox{{\tt
		\quad\quad  {\small \rm #1}}}%
	\theight=\ht0
	\advance\theight by \lineskip
	\kern -\theight \vbox to
	\theight{\rightline{\rlap{\box0}}%
	\vss}%
	}}%

\ShowLabelsfalse

\renewcommand{\sec}[2]{\section{#2}\label{S:#1}%
	\ifShowLabels \TeXref{{S:#1}} \fi}
\newcommand{\ssec}[2]{\subsection{#2}\label{SS:#1}%
	\ifShowLabels \TeXref{{SS:#1}} \fi}

\newcommand{\refs}[1]{Section ~\ref{S:#1}}
\newcommand{\refss}[1]{Section ~\ref{SS:#1}}
\newcommand{\reft}[1]{Theorem ~\ref{T:#1}}
\newcommand{\refl}[1]{Lemma ~\ref{L:#1}}
\newcommand{\refp}[1]{Proposition ~\ref{P:#1}}
\newcommand{\refc}[1]{Corollary ~\ref{C:#1}}
\newcommand{\refd}[1]{Definition ~\ref{D:#1}}
\newcommand{\refr}[1]{Remark ~\ref{R:#1}}
\newcommand{\refe}[1]{\eqref{E:#1}}

\newenvironment{thm}[1]%
	{ \begin{Thm} \label{T:#1}  \ifShowLabels \TeXref{T:#1} \fi }%
	{ \end{Thm} }

\newcommand{\th}[1]{\begin{thm}{#1} }
\renewcommand{\eth}{\end{thm} }

\newenvironment{lemma}[1]%
	{ \begin{Lem} \label{L:#1}  \ifShowLabels \TeXref{L:#1} \fi }%
	{ \end{Lem} }
\newcommand{\lem}[1]{\begin{lemma}{#1}}
\newcommand{\elem}{\end{lemma}}

\newenvironment{propos}[1]%
	{ \begin{Prop} \label{P:#1}  \ifShowLabels \TeXref{P:#1} \fi }%
	{ \end{Prop} }
\newcommand{\prop}[1]{\begin{propos}{#1}}
\newcommand{\eprop}{\end{propos}}

\newenvironment{corol}[1]%
	{ \begin{Cor} \label{C:#1}  \ifShowLabels \TeXref{C:#1} \fi }%
	{ \end{Cor} }
\newcommand{\cor}[1]{\begin{corol}{#1}}
\newcommand{\ecor}{\end{corol}}

\newenvironment{defeni}[1]%
	{ \begin{Def} \label{D:#1}  \ifShowLabels \TeXref{D:#1} \fi }%
	{ \end{Def} }
\newcommand{\defe}[1]{\begin{defeni}{#1}}
\newcommand{\edefe}{\end{defeni}}

\newenvironment{remark}[1]%
	{ \begin{Rem} \label{R:#1}  \ifShowLabels \TeXref{R:#1} \fi }%
	{ \end{Rem} }
\newcommand{\rem}[1]{\begin{remark}{#1}}
\newcommand{\erem}{\end{remark}}

\newcommand{\eq}[1]%
	{ \ifShowLabels \TeXref{E:#1} \fi
	   \begin{equation} \label{E:#1} }
\newcommand{\eeq}{ \end{equation} }

\newcommand{\prf}{ \begin{pf} }
\newcommand{\epr}{ \end{pf} }

\newcommand\alp{\alpha}		
\newcommand\bet{\beta}
		
\newcommand\del{\delta}		\newcommand\Del{\Delta}
\newcommand\eps{\varepsilon}

\newcommand\tet{\theta}		

\newcommand\kap{\kappa}
\newcommand\lam{\lambda}		\newcommand\Lam{\Lambda}
\newcommand\sig{\sigma}		
\newcommand\vphi{\varphi}
\newcommand\ome{\omega}		\newcommand\Ome{\Omega}

\newcommand\calF{{\cal{F}}}

\newcommand\calH{{\cal{H}}}

\newcommand\calL{{\cal{L}}}
\newcommand\calM{{\cal{M}}}
\newcommand\calN{{\cal{N}}}

\newcommand\calP{{\cal{P}}}
\newcommand\calQ{{\cal{Q}}}
\newcommand\calR{{\cal{R}}}

\newcommand\RR{\Bbb{R}}

\newcommand\ZZ{\Bbb{Z}}

\newcommand\CC{\Bbb{C}}

\newcommand\nek{,\ldots,}
\newcommand\sdp{\times \hskip -0.3em {\raise 0.3ex
\hbox{$\scriptscriptstyle |$}}} 

\newcommand\Dom{\operatorname{Dom}}

\newcommand\GL{\operatorname{GL}}

\newcommand\IM{\operatorname{Im}}

\newcommand\ind{\operatorname{ind}}

\newcommand\Ker{\operatorname{Ker}}

\newcommand\rk{\operatorname{rk}}

\newcommand\supp{\operatorname{supp}}

\newcommand\oJ{{\overline{J}}}

\newcommand\tilH{{\widetilde{H}}}

\newcommand\tilOme{{\widetilde{\Ome}}}

\theoremstyle{plain}

\newcommand{\F}{\calF}
\renewcommand{\d}{\text{\( \partial\)}}
\newcommand{\w}{\text{\( \ome\)}}
\newcommand{\we}[1]{\text{\( \Ome_{(2)}^{#1}(E,\F)\)}}
\newcommand{\wess}[1]{\text{\( \tilOme_{(2)}^{#1}(E,\F;\n_{1,\alp})\)}}
\newcommand{\wes}[1]{\text{\( \Ome^{#1}(E,\F)\)}}

\renewcommand{\t}{\text{\( \tau_{\sqrt{t}}\)}}
\newcommand{\r}{\text{\( r^*_{\sqrt{t}}\)}}

\newcommand{\nd}{\text{\( \Del^C_{t,\alp}\)}}
\newcommand{\ndel}[2]{\text{\( \Del^{C,#1}_{#2,\alp}\)}}
\newcommand{\ndp}{\text{\( \Del^{C,p}_{1,\alp}\)}}

\renewcommand{\H}{\text{\( \Del^p_{t,\alp}\)}}

\newcommand{\n}{\nabla}

\renewcommand{\b}{\bullet}

\newcommand{\m}{m}

\begin{document}
\title[Novikov type inequalities]
	{Novikov type inequalities for differential forms
		with non-isolated zeros}
\author{Maxim Braverman} \author{Michael Farber}
\address{School of Mathematical Sciences\\
Tel-Aviv University\\
Ramat-Aviv 69978, Israel}
\email{maxim@@math.tau.ac.il, \,farber@@math.tau.ac.il}
\thanks{The research was supported by grant No. 449/94-1 from the
Israel Academy of Sciences and Humanities.}

\maketitle

\begin{abstract}
We generalize the Novikov inequalities for 1-forms in two different
directions:
first, we allow non-isolated critical points (assuming that they are
non-degenerate in the sense of R.Bott), and, secondly, we strengthen
the inequalities by means of twisting
by an arbitrary flat bundle.
The proof uses Bismut's modification of the Witten deformation of the de
Rham complex;
it is based on an explicit estimate
on the lower part of the spectrum of the corresponding Laplacian.

In particular, we obtain a new
analytic proof of the degenerate Morse inequalities of Bott.
\end{abstract}

\sec{introd}{Introduction}

\ssec{novin}{The Novikov inequalities} Let $M$ be a closed manifold of
dimension $n$.
In \cite{n1,n2}, S.P.Novikov associated to any real cohomology class
$\xi\in H^1(M,\RR)$ a sequence of numbers $\beta_0(\xi),\dots ,\beta_n(\xi)$
and proved that for any closed 1-form $\omega$ on $M$, having non-degenerate
critical points, the following inequalities hold
\eq{novin}
	m_p(\omega)\ge \beta_p(\xi), \qquad p=0,1,2,\dots
\end{equation}
Here $\xi=[\omega]\ \in H^1(M,\RR)$ is the cohomology class of $\omega$ and
$m_p(\omega)$ denotes the number of critical points of $\omega$ having
Morse index $p$.
Note also, that there are slightly stronger inequalities
\eq{novin'}
	\sum_{i=0}^{p}(-1)^im_{p-i}(\omega)\
\ge\  \sum_{i=0}^{p}(-1)^i\beta_{p-i}(\xi), \qquad p=0,1,2,\dots,
\end{equation}
cf. \cite{far}.
In the case when the form $\omega$ is {\it exact} (i.e.
$\omega=df$ where $f$ is a non-degenerate Morse function) then $\xi = 0$
and the Novikov inequalities \refe{novin'} turn into the classical Morse
inequalities.

A recent survey of the theory of Novikov inequalities may be found in
\cite{p2}.

\ssec{main}{Formulation of the main result} In this paper we generalize
the Novikov inequalities for 1-forms in two different
directions:
firstly, we allow non-isolated critical points and, secondly, we strengthen
the inequalities by means of twisting
by an arbitrary flat vector bundle.

Let $M$ be a closed smooth manifold with a fixed flat complex vector bundle
$\F$. Let $\omega$ be a smooth closed real valued 1-form on $M$, $d\omega=0$,
which is assumed to be {\it non-degenerate in the sense of R.Bott}
\cite{bott}. This means that the points of $M$, where the form
$\omega$ vanishes form a submanifold of $M$ (called the {\em
critical points set} $C$ of $\omega$) and the {\it Hessian} of $\omega$ is
{\it non-degenerate on the normal bundle to} $C$.

In order to make clear this definition, note that if we fix a tubular
neighborhood $N$ of $C$ in $M$, then the monodromy of $\omega$ along
any loop in $N$ is obviously zero. Thus there exists a unique real
valued smooth function $f$ on $N$ such that $df=\omega_{|_N}$ and
$f_{|_C}=0$. The {\it Hessian} of $\omega$ is then defined as the
Hessian of $f$.

Let $\nu(C)$ denote the normal bundle of $C$ in $M$. Note that
$\nu(C)$ may have different dimension over different connected
components of $C$.  Since Hessian of $\omega$ is non-degenerate, the
bundle $\nu(C)$ splits into the Whitney sum of two subbundles
\eq{n=n+n-}
	\nu(C)\ =\ \nu^+(C)\ \oplus \nu^-(C),
\end{equation}
such that the Hessian is strictly positive on $\nu^+(C)$ and strictly
negative on $\nu^-(C)$. Here again, the dimension of the bundles
$\nu^+(C)$ and $\nu^-(C)$ over different connected components of the
critical point set may be different.

For every connected component $Z$ of the critical point set $C$, the
dimension of the bundle $\nu^-(C)$ over $Z$ is called the {\em index}
of $Z$ (as a critical submanifold of \w) and is denoted by $\ind(Z)$.
Let $o(Z)$ denote
the {\it orientation bundle of $\nu^-(C)_{|_Z}$, considered as a flat
line bundle}. Consider the {\it twisted Poincar\'e polynomial} of component
$Z$
\eq{P-pol}
	\calP_{Z,\F}(\lambda)\ =\
	   \sum_{i=0}^{\dim Z} \lambda^i\dim_{\CC} H^i(Z,\F_{|_Z}\otimes o(Z))
\end{equation}
(here $H^i(Z,\F_{|_Z}\otimes o(Z))$ denote the cohomology of $Z$ with
coefficients in the flat vector bundle $\F_{|_Z}\otimes o(Z)$)
and define using it the following {\it Morse counting polynomial}
\eq{M-pol}
	\calM_{\omega,\F}(\lambda)\ =
		\ \sum_{Z} \lambda^{\ind(Z)} \calP_{Z,\F}(\lambda),
\end{equation}
where the sum is taken over all connected components $Z$ of $C$.

On the other hand, with one-dimensional cohomology class
$\xi=[\omega]\in \-H^1(M,\RR)$ and the flat vector bundle $\F$, one can
associate canonically
the {\em Novikov counting polynomial}
\eq{N-pol}
	\calN_{\xi,\F}(\lambda)\ =\ \sum_{i=0}^n \lambda^i \beta_i(\xi,\F),
\end{equation}
where $\beta_i(\xi,\F)$ are generalizations of the Novikov numbers,
cf. \refd{novnum}. Note, that if $\xi=0$ and $\F$ is the trivial line
bundle, then $\calN_{\xi,\F}(\lam)$ coincides with the Poincar\'e polynomial
of $M$.

The following is our principal result.
\th{main} In the situation described above, there exists a polynomial
  $$
	\calQ(\lambda)= q_0 +  q_1\lambda + q_2\lambda^2 + \dots
  $$
  with non-negative
  integer coefficients $q_i\ge 0$, such that
  \eq{main}
	\calM_{\omega, \F}(\lambda)\ -\ \calN_{\xi,\F}(\lambda)\ =\
			(1+\lambda)\calQ(\lambda).
  \end{equation}
\eth

The main novelty in this theorem is that it is applicable to the case
of 1-forms with non-isolated singular points. Thus, we obtain, in
particular, a new proof of the degenerate Morse inequalities of
R.Bott.  Moreover, \reft{main} provides a generalization of the
Morse-Bott inequalities to the case of an arbitrary flat vector bundle
$\F$; this generally produces stronger inequalities  as shown in
\refss{example}.

Next, we are going to point out the following corollary.

\begin{Cor}[Euler-Poincar\'e theorem] \label{C:Eu-Poin}
  Under the conditions of \reft{main}, the Euler
  characteristic of $M$ can be computed as
  \eq{Eu-Poin}
	\chi(M)\ =\ \sum_{Z} \ (-1)^{\ind(Z)}\chi(Z),
  \end{equation}
  where the sum is taken over all connected components $Z\subset C$.
\end{Cor}
The corollary is obtained from \refe{main} by substituting $\lam=-1$
and observing that
\eq{Eu-Poin'}
	\calM_{\w,\F}(-1)=d \cdot \sum_{Z}\  (-1)^{\ind(Z)}\chi(Z),
	\qquad \calN_{\xi,\F}(-1)=d \cdot \chi(M),
\end{equation}
where $d=\dim\F$; the last equality in \refe{Eu-Poin'}  follows
immediately from \refd{novnum}.

\ssec{isol}{The case of isolated critical points}
Let's consider the special case when all critical points of $\omega$
are isolated.  Then the Morse counting polynomial \refe{M-pol} takes
the form
\eq{M-pol'}
	\calM_{\omega,\F}(\lambda)\ =\ d\cdot \sum_{p=0}^n \lam^p m_p(\omega),
\end{equation}
where $d=\dim \F$ and $m_p(\omega)$ denotes the number of critical points
of $\omega$ of index $p$.  \reft{main} gives in this case the
inequalities
\eq{novin''}
	\sum_{i=0}^{p}(-1)^im_{p-i}(\omega)\
		\ge\  d^{-1}\cdot \sum_{i=0}^{p}(-1)^i\beta_{p-i}(\xi, \F),
				\qquad p=0,1,2,\dots
\end{equation}
The last inequalities coincide with the Novikov inequalities
\refe{novin'} in the special case when $\F=\RR$ with the trivial flat
structure. Easy examples described in \refss{example}, show that using of
the flat vector bundle $\F$ gives sharper estimates in general, than
the standard approach with $\F=\RR$.

On the other hand, \refe{novin''} clearly generalizes the Morse type
inequalities obtained by S.P.Novikov in \cite{n3}, using Bloch
homology (which correspond to the case, when $[\omega]=0\in
H^1(M,\RR)$ in \refe{novin''}).

\ssec{method}{The method of the proof}
Our proof of \reft{main} is based on a slight modification of the Witten
deformation  \cite{wi}  suggested by Bismut \cite{bis}
in his proof of the degenerate Morse inequalities of Bott.
However  our proof is rather different from \cite{bis} even in the case
$[\w]=0$. We entirely avoid the probabilistic
analysis of the heat kernels, which is the most difficult part of \cite{bis}.
Instead, we give an explicit estimate on
the number of the "small" eigenvalues of the deformed Laplacian.
We now will explain briefly the main steps of the proof.

Let $\Ome^\b(M,\F)$ denote the space of smooth differential forms on
$M$ with values in $\F$. In \refs{manif}, we introduce a 2-parameter
deformation
\eq{n-talp}
	\n_{t,\alp}:\Ome^\b(M,\F)\to \Ome^{\b+1}(M,\F),
		\qquad t,\alp \in \RR
\end{equation}
of the covariant derivative $\n$, such that, for large values of
$t,\alp$ the Betti numbers of the deformed de Rham complex
$\big(\Ome^\b(M,\F),\n_{t,\alp}\big)$ are equal to the Novikov numbers
$\bet_p(\xi,\F)$.  Outside of a small tubular neighborhood of the
critical points set $C$ of \w\ the differential \refe{n-talp} is given
by the formula
$$
	\n_{t,\alp}:\tet\mapsto \tet+t\alp\w\wedge\tet,      \qquad
						\tet\in\Ome^\b(M,\F).
$$
Next we construct a special Riemannian metric $g^M$ on $M$. In fact,
we, first, chose a Riemannian metric on the normal bundle $\nu(C)$
and, then, extend it to a metric on $M$. We also choose a Hermitian
metric $h^\F$ on $\F$.  Let us denote by $\Del_{t,\alp}$ the Laplacian
associated with the differential \refe{n-talp} and with the metrics
$g^M,h^\F$.

Fix $\alp>0$ sufficiently large. It turns out that, when $t\to\infty$,
the eigenfunctions of $\Del_{t,\alp}$ corresponding to ``small''
eigenvalues localize near the critical points set $C$ of \w. Hence,
the number of the "small" eigenvalues of $\Del_{t,\alp}$ may be
calculated by means of the restriction of $\Del_{t,\alp}$ on a tubular
neighborhood of $C$. This neighborhood may be identified with a
neighborhood of the zero section of the normal bundle $\nu(C)$ to
$C$. We are led, thus, to study of a certain Laplacian on
$\nu(C)$. The latter Laplacian may be decomposed as
$\bigoplus_Z\Del^Z_{t,\alp}$ where the sum ranges over all connected
components of $C$ and $\Del^Z_{t,\alp}$ is a Laplacian on the normal
bundle $\nu(Z)=\nu(C)_{|_Z}$ to $Z$. We denote by $\Del^{Z,p}_{t,\alp}
\ (p=0,1,2,\dots)$ the restriction of $\Del^Z_{t,\alp}$ on the space
of $p$-forms.

The operator \refe{n-talp} is constructed so that the spectrum of
$\Del^Z_{t,\alp}$ does not depend on $t$. Moreover, if $\alp>0$ is
sufficiently large, then
\eq{bismut'}
	\dim\Ker\Del^{Z,p}_{t,\alp}=\dim H^{p-\ind(Z)}(Z,\F_{|_Z}).
\end{equation}
In the case when $\F$ is a trivial line bundle, the equation
\refe{bismut'} is proven by Bismut \cite[Theorem 2.13]{bis}. We prove
\refe{bismut'} in \refs{L2cohom}.  Note that our proof is rather
different from \cite{bis}.

Let $E^\b_{t,\alp} \ (p=0,1\nek n)$ be the subspace of $\Ome^\b(M,\F)$
spanned by the eigenvectors of $\Del_{t,\alp}$ corresponding to the
``small'' eigenvalues. The cohomology of the deformed de Rham complex
$\big(\Ome^\b(M,\F),\n_{t,\alp}\big)$ may be calculated as the
cohomology of the subcomplex $\big(E^\b_{t,\alp},\n_{t,\alp}\big)$.

We prove (\reft{spectr}) that, if the parameters $t$ and $\alp$ are
large enough, then
\eq{k=h'}
	\dim E^p_{t,\alp}= \sum_{Z}\dim\Ker\Del^{Z,p}_{t,\alp},
\end{equation}
where the sum ranges over all connected components $Z$ of $C$. The
\reft{main} follows now from \refe{bismut'},\refe{k=h'} by standard
arguments (cf. \cite{bott2}).
\rem{hel-sj}In \cite{hs3}, Helffer and Sj\"ostrand  gave a
  very elegant analytic proof of the degenerate Morse inequalities of
  Bott. Though they also used the ideas of \cite{wi}, their method is
  completely different from \cite{bis}.  It is not clear if this method may
  be applied to the case $\xi\not=0$.
\erem

\ssec{contents}{Contents} The paper is organized as follows.

In \refs{nov}, we define a slightly generalized version of the Novikov
numbers associated to a cohomology class $\xi\in H^1(M,\RR)$ and a
flat vector bundle $\F$. Here we also discuss some examples.

In \refs{bundle}, we recall the construction of Bismut's deformation
of the de Rham complex on a fiber bundle and discuss the  spectral
properties of the corresponding Laplacian. At the end of the section we
state \reft{bismut} which calculates the kernel of
this Laplacian.

In \refs{L2cohom}, we  prove \reft{bismut}.

In \refs{manif}, we define the Bismut deformation of the Laplacian on $M$,
taking into account that \w\ is not cohomologious to 0. Then we
prove that the number of the   eigenvalues of this Laplacian which tend to 0
as the parameter tends to infinity is equal to the dimension of the kernel of
some Laplacian on $\nu(Z)$.

In \refs{proof}, we prove the main \reft{main}.

\

The results contained in this paper were announced in \cite{bf}.


\sec{nov}{The Novikov numbers}

In this section we recall the definition and the main properties of the
Novikov numbers \cite{n1,n2} associated to a cohomology class
$\xi\in H^1(M,\RR)$. In fact, we define these numbers in a slightly more
general situation. Our point of view is motivated by the study of
deformations of elliptic complexes in \cite{far1}. Roughly speaking,
any one-dimensional cohomology class defines an analytic deformation of
the twisted de Rham complex and the Novikov numbers are the instances of
the natural invariants of such deformations, cf. \cite{far1}.


\ssec{1.1}{The Novikov deformation} Let $M$ be a closed smooth
manifold and let $\F$ be a complex flat vector bundle over $M$. We
will denote by
$$
	\nabla:\Omega^\b(M,\F)\to \Omega^{\b+1}(M,\F)
$$
the covariant derivative on $\F$.

Given a closed 1-form $\omega\in \Omega^1(M)$ on $M$ with real values,
it determines a family of connections on $\F$ (the {\em Novikov deformation})
parameterized by the real numbers $t\in \RR$
\eq{nab-t}
	\nabla_t:\Omega^i(M,\F)\to \Omega^{i+1}(M,\F),
\end{equation}
where
\eq{nub-t'}
	\nabla_t\theta=\nabla\theta+t\omega\wedge\theta,\qquad
		\theta\in\Omega^\bullet(M,\F).
\end{equation}
All the connections $\nabla_t$ are flat, i.e. $\nabla_t^2=0$,
if the form $\omega$ is closed, $d\omega=0$.

We can view the obtained complex as follows. For $t\in\RR$,
let $\rho_t$ denote the flat real line bundle over $M$ with the monodromy
representation $\rho_t:\pi_1(M)\to \RR^\ast$ given by the formula
\eq{novrep}
	\rho_t(\gamma)=
		\exp(-t\int_\gamma\omega)\ \ \in \ \RR^\ast,\qquad \gamma\in
							\pi_1(M).
\end{equation}
Then $\nabla_t$ can be considered as the covariant derivative
on the flat bundle $\F\otimes\rho_t$.

Note that changing $\omega$ by a cohomologious 1-form determines a gauge
equivalent connection $\nabla_t$ and so the cohomology
$H^\b(M,\F\otimes\rho_t)$ depends only on the cohomology class
$\xi= [\omega]\in H^1(M,\RR)$ of $\omega$.

The dimension of the cohomology $H^i(M,\F\otimes\rho_t)$ is an integer
valued function of $t\in\RR$ having the following behavior. There exists a
discrete subset $S\subset\RR$ (i.e. each of its points is isolated)
such that the dimension
$\dim \ H^i(M,\F\otimes\rho_t)$ is {\it constant for $t\notin S$}
(the corresponding value of the dimension we will call the
{\it background value}) and for $t\in S$ the dimension of
$H^i(M,\F\otimes\rho_t)$ is
{\it greater} than the background value.
Cf., for example, \cite[Theorem 2.8]{far1},
where a more precise information for the case of
elliptic complexes is given.
The subset $S$ above will be called
{\it the set of jump points}.

\defe{novnum} For each $i=0,1\nek n$, the background value of the
  dimension of $H^i(M,\F\otimes\rho_t)$ is called the $i$-th
  {\em Novikov number} $\beta_i(\xi,\F)$.
\edefe

The novelty here is in introduction of the flat vector bundle $\F$;
the standard definition uses the trivial line bundle (over $\CC$)
instead of our $\F$.  The importance of this generalization
will be explained in \refss{example} below.

\lem{1.2} The set of jump points $S$ is finite.
\elem
\prf
Given class $\xi\in H^1(M,\RR)$, consider the set $\cal R(\xi)$ of all
cohomology classes $\rho\in H^1(M,\RR)$ such that the the
corresponding period map $\rho_\ast: H_1(M,\ZZ)\to \RR$ vanishes on
the kernel of the period map $\xi_\ast: H_1(M,\ZZ)\to \RR$.  The image
of the period map $\xi_\ast: H_1(M,\ZZ)\to \RR$ is isomorphic to
$\ZZ^l$ for some $l$ (which is called the {\em degree of irrationality
of $\xi$}). Any element $\rho\in \cal R(\xi)$ determines a
representation of the fundamental group
$$
	\pi=\pi_1(M) \ \to \RR^{>0},\qquad\text{where}\quad
			g\mapsto \exp(-\rho_\ast[g])\ \in\RR^{>0}.
$$
This allows to identify $\cal R(\xi)$ with $(\RR^\ast)^l$ and thus we
introduce the structure of affine algebraic variety in $\cal R(\xi)$.

Consider the function $\dim_{\CC} H^i(M,\F\otimes\rho)$ as a function of
$\rho\in \cal R(\xi)$. The standard arguments show that there exist an
algebraic  $V\subsetneq (\RR^\ast)^l$ such that
$\dim_{\CC} H^i(M,\F\otimes\rho)$ is constant for all $\rho\notin V$
(cf. \cite[Ch.3 \S 12]{harts}).

If we identify $\cal R(\xi)$ with $(\RR^\ast)^l$ as explained above,
then the point $\xi$ will be  represented by a vector
$(e^{a_1},e^{a_2},\dots,e^{a_l})$ with the real numbers $a_1, a_2,\dots,a_l$
linearly independent over $\Bbb Q$. Then, for $t\in\RR$, the class
$t\xi$ (describing the 1-dimensional local system $\rho_t$) is
represented by the vector $(e^{ta_1},e^{ta_2},\dots,e^{ta_l})\in (\RR^\ast)^l$.
It is easy to see that the curve
$$\psi(t)\ =\ (e^{ta_1},e^{ta_2},\dots,e^{ta_l})\in (\RR^\ast)^l$$
has the following property: for any algebraic subvariety
$W\subset (\RR^\ast)^l$ there exists a constant $C>0$ such that
$\psi(t)\notin W$ for $|t|>C$.

This show that $\dim_{\CC} H^i(M,\F\otimes\rho_t)$ is equal to the background
value for $t$ sufficiently large. Thus the set of jump points $S$ is finite.
\epr

\ssec{1.3}{The monodromy of the deformed connection $\nabla_t$}
Let $\F_0$ denote the fiber of $\F$ over the base point $\ast$ of
$M$. Let $\phi:\pi_1(M)\to \GL_{\CC}(\F_0)$ denote the monodromy
representation of the flat bundle $\F$. Then the monodromy
representation of the flat vector bundle $\F\otimes \rho_t$ is given
by
$$
	g\mapsto \exp(-t\int_g\omega) \phi(g)\ \in \GL_{\CC}(\F_0)
				\qquad\text{for}\qquad g\in\pi_1(M).
$$
This formula describes explicitly the {\it deformation of the
monodromy representation}.

\ssec{1.4}{Computation of the Novikov numbers in terms of the spectral
sequence}
One may compute the Novikov numbers $\beta_i(\xi,\F)$ by means of the
cell structure of $M$ as the dimension of the homology of a local
system over $M$ determined by the deformation \refe{nub-t'}; the
dimension here is understood over the field of germs of meromorphic
curves in ${\CC}$. This was explained in \cite{far1} and in \cite{fl}
in terms of  the {\em germ complex of the deformation}; we will not
repeat this construction here. It leads naturally to the computation
of the Novikov numbers by means of a spectral sequence, cf. Theorems
2.8 and 6.1 in \cite{far1}. It is important to point out that both
these approaches (the one, based on the germ complex, and the spectral
sequence) are able to find the Novikov numbers {\em starting from an
arbitrary value $t=t_0$ and using infinitesimal information on the
deformation}.

Note, that a spectral sequence of a similar nature appears in \cite{p2}.

We will briefly describe the spectral sequence.
Fix an arbitrary value $t=t_0\in\RR$. Then there exists a spectral sequence
$E^\ast_r$, $\ r\ge 1$ with the following properties:
\begin{enumerate}
\item The initial term of the spectral sequence coincides with the cohomology
at the chosen point $t=t_0$:
$$E^i_1\ =H^i(M,\F\otimes\rho_{t_0});$$
\item For large $r$ all differentials of the spectral sequence
$d_r: E^i_r\to E^{i+1}_r$ vanish and the limit term
$E^i_\infty$ is isomorphic to the background cohomology;
\item The first differential
$$H^i(M,\F\otimes\rho_{t_0})\ \to H^{i+1}(M,\F\otimes\rho_{t_0})$$
is given by multiplication by $\xi\in H^1(M,\RR)$;
\item The higher differentials are given by the iterated Massey products
with $\xi$.
\end{enumerate}

The spectral sequence is constructed as follows, cf.\cite[Section 6.1]{far1}.
Denote by $Z^i_r$ the set of polynomials of the form
$$
	f(t)=f_0+tf_1+t^2f_2+\cdots +t^{r-1}f_{r-1}
$$
where $f_j\in \Omega^i(M,\F\otimes\rho_{t_0})$ such that $\nabla_tf(t)$
is divisible by $(t-t_0)^r$. Then
$$
	E^i_r\ =\ Z^i_r/(tZ^i_{r-1}+t^{1-r}\nabla_tZ^{i-1}_{r-1})
$$
and the differential
$$
	d_r: E^i_r\ \to\ E^{i+1}_r
$$
is induced by the action of $t^{-r}\nabla_t$. See \cite[Section 6]{far1},
for more detail.

\cor{1.5} If $\xi\in H^1(M,\RR)$ is a non-zero class then
  the zero-dimensional Novikov number $\beta_0(\xi,\F)$ vanishes.
\ecor

It follows directly by applying the spectral sequence.

\ssec{example}{Some examples} Here we will produce examples, where
the Novikov
numbers twisted by a flat vector bundle $\F$ (as defined above) give greater
values (and thus stronger inequalities) than the usual Novikov numbers
(where $\F=\RR \text{ or } {\CC}$, cf. \cite{n1,n2,p2}).

Let $k\subset S^3$ be a smooth knot and let the 3-manifold $X$ be the result
of $1/0$-surgery on $S^3$ along $k$. Note that the one-dimensional homology
group of $X$ is infinite cyclic and thus for any complex number
$\eta\in{\CC}$, $\eta\ne 0$, there is a complex flat line bundle over
$X$ such that the monodromy with
respect to the generator of $H_1(X)$ is $\eta$. Denote such flat bundle by
$\F_\eta$.

Note also, that for $\eta\neq 1$, the dimension of $H^1(X,\F_\eta)$ is
equal to the multiplicity of $\eta$ as a root of the Alexander
polynomial of the knot $k$, cf. \cite{rol}.  Thus, by a choice of the
knot $k$ and the number $\eta\in{\CC}^\ast$, we may make the group
$H^1(X,\F_\eta)$ arbitrarily large, while $H^1(X,{\CC})$ is always
one-dimensional. (One may, for example, take multiple connected sum of
many copies of a knot with non-trivial Alexander polynomial).

Consider now the 3-manifold $M$ which is the connected sum
$$
	M\ =\ X\ \#\ (S^1\times S^2).
$$
Thus $M=X_+\cup X_-$ where $X_+\cap X_-=X_0=S^2$ and
$$
	X_+=X-\text\{disk\}\quad\text{and}\quad
			X_-=(S^1\times S^2)-\text\{disk\}.
$$
Consider a flat complex line bundle $\F$ over $M$ such that its restriction
over $X_+$ is isomorphic to $\F_\eta|_{X_+}$. Consider the class
$\xi\in H^1(X,\RR)$
such that its restriction onto $X_+$ is trivial and its restriction to
$X_-$ is the generator.

We want to compute the Novikov number $\beta_1(\xi,\F)$. By using the
Mayer-Vietoris sequence, we obtain, for generic $t$,
$$
	H^1(M,\F\otimes \rho_t)\ \simeq\ H^1(X_+,\F_\eta)\oplus
			H^1(X_-,\F\otimes\rho_t)
$$
For generic $t$, the second term vanishes and thus we obtain
$$
	\beta_1(\xi,\F)\ =\ \dim_{{\CC}} H^1(X,\F_\eta)
$$
As we noticed above, this number can be arbitrarily large, while
$\dim_{\CC} H^1(M,{\CC})=2$.

\sec{bundle}{The Bismut Laplacian on a fiber bundle}

In this section we describe a version of the Bismut deformation of the
Laplacian on the bundle $E$, considered as a non-compact manifold. In
\refs{manif} we will apply this construction to the case when $E$ is
the normal bundle to a connected submanifold of the set of critical
points of \w.

\ssec{data}{Description of the data} Let $E=E^+\oplus E^-$ be a
$\ZZ_2$-graded finite dimensional vector bundle over a compact
connected manifold $Z$.

Suppose that $\F$ is a flat vector bundle over $E$. We denote by \wes{\b}
the space of differential forms on $E$ with values in $\F$ and by
$\n:\wes{\b}\to \wes{\b+1}$ the corresponding covariant derivative operator.


\ssec{split}{A splitting of the tangent space} Fix an Euclidean metric
$h^E$ on the bundle $E$ (i.e. a smooth fiberwise metric) such that
$E^+$ and $E^-$ are orthogonal.  For a vector $y\in E$, we will
denote by $|y|$ its norm with respect to the metric $h^E$.

Next choose an Euclidean connection $\n^E$ on $E$ which preserves the
decomposition $E=E^+\oplus E^-$. Then the tangent space $TE$ splits
naturally into
\eq{split}
	TE=T^HE\oplus T^VE,
\end{equation}
where $T^VE$ is the set of the vectors  in $TE$ which are tangent to the
fibers of $E$ (the {\em vertical vectors}), and $T^HE$ is the set
of {\em horizontal vectors} in $TE$.


\ssec{metric}{A Riemannian metric on $E$}In this section we consider
$E$ as a non-compact manifold. Our aim is to introduce a Riemannian
metric on $E$.  Let $\pi$ be the projection $E\to Z$.  If $y\in E$,
then $\pi_*$ identifies $T_y^HE$ with $T_{\pi(y)}Z$.  Choose any
Riemannian metric $g^Z$ on $Z$. Then $T_y^HE$ is naturally endowed
with the metric $\pi^*g^Z$. Also $T^VE$ and $E$ can be naturally
identified. Hence, the metric $h^E$ on $E$ induces a metric on
$T^VE$. We still denote this metric by $h^E$ and we define the
Riemannian metric
\eq{gE}
	g^E=h^E\oplus \pi^*g^Z
\end{equation}
on $TE$, which coincides with $h^E$ on $T^VE$, with $\pi^*g^Z$
on $T^HV$ and such that $T^HE$ and $T^VE$ are orthogonal.
Note that the metric $g^E$ depends upon the choices of $h^E,g^Z$ and
$\n^E$.

\ssec{metF}{An Euclidean metric on $\calF$}
We will identify the manifold $Z$ with the zero section of $E$. Let
$\F_{|_Z}$ denote the restriction of $\F$ on $Z$. Fix an arbitrary
Euclidean metric $h$ on $\F_{|_Z}$. The flat connection on $\F$
defines a trivialization of $\F$ along the fibers of $E$ and,
hence, gives a natural extension of $h$ to an Euclidean metric $h^\F$
on $\F$ which is flat along the fibers of $E$.

%

\ssec{bigrad}{A bigrading on the space of differential forms} The
metrics $g^E, h^\F$ define an $L_2$-scalar product on the space of
differential forms on $E$ with values in $\F$. Let $\we{}=\bigoplus
\we{p}$ denote the Hilbert space of square integrable differential forms on
$E$ with values in $\F$.

The splitting \refe{split} and the corresponding decomposition
$T^*E=(T^HE)^*\oplus (T^VE)^*$ induce a bigrading on $\we{}$ by
\eq{bigr}
	\we{}=\bigoplus_{i,j} \we{i,j},
\end{equation}
where \we{i,j}\ is the space of square integrable sections of
$$
	\Lam^i\big((T^HE)^*\big)\otimes\Lam^j\big((T^VE)^*\big)\otimes \F.
$$
Here, as usual, we denote by $\Lam(V)=\bigoplus \Lam^p(V)$ the
exterior algebra of a vector space $V$.

For $s\not= 0$, let $\tau_s$ be the map from \we{} onto itself, which
sends $\alp\in \we{i,j}$ to $s^j\alp$. Of course, $\tau_s$ depends
upon the connection on $E$ given by  \refe{split}.

\ssec{bismut}{The Bismut complex on a fiber bundle}
In this section we produce a family of complexes depending on $t>0$,
which is a special case of the construction of Bismut
\cite[Section 2(b)]{bis}. Note that our notations are slightly
different from \cite{bis}.

Let $f:E\to \RR$ denote the function such that its value on a vector
$y=(y^+,y^-)\in E^+\oplus E^-$ is given by
\eq{fE}
	f(y)=\frac{|y^+|^2}2-\frac{|y^-|^2}2.
\end{equation}

Recall that $\n:\we{\b}\to \we{\b+1}$ denotes the covariant
differential operator determined by the flat connection on $\F$.
Following Bismut \cite{bis}, we define a family of differentials
\eq{dE}
	\n_{t,\alp}=(\t)^{-1}e^{-\alp tf}\n e^{\alp tf}\t,
			\qquad	t>0,\ \alp>0.
\end{equation}

\ssec{dE'}{An alternative description of $\n_{t,\alp}$}
We will also  need another description of the differential \refe{dE}
(cf. \cite[Remark 2]{bis}).

For $s>0$, let $r_s:E\to E$ be the multiplication by $s$, i.e
$r_sy=sy$ for any $y\in E$. Recall that, if $y\in E$, we denote by
$\F_y$ the fiber of $\F$ over $y$. The flat connection on $\F$ gives a
natural  identification of  the fibers $\F_y$ and $\F_{sy}$.
Hence, the map $r_s:E\to E$ defines the ``pull-back'' map
$$
	r^*_s:\we{\b}\to \we{\b}.
$$
Note that $r^*_s$ preserves the bigrading \refe{bigr}
and, hence, commutes with \t.
\lem{de-de'} For any $t>0, \alp >0$
  \eq{dE'}
	\n_{t,\alp}=(\t)^{-1}\r e^{-\alp f}\n e^{\alp f}(\r)^{-1}\t.
  \end{equation}
\elem
\prf
  Since \r\ commutes with $\n$ and with \t, it is enough to show that
  \eq{tft}
	\r\circ f\circ (\r)^{-1}=tf,
  \end{equation}
  where $f$ is identified with the operator of multiplication by $f$.
  The equality \refe{tft} follows immediately from \refe{fE}.
\epr

\defe{lapl}The {\it Bismut Laplacian} $\Del_{t,\alp}$ of the bundle
$E$ associated to the metrics $g^E, h^{\F}$ is defined by the formula
\eq{lapl}
    \Del_{t,\alp}=\n_{t,\alp}\n^*_{t,\alp}+\n^*_{t,\alp}\n_{t,\alp},
\end{equation}
where $\n^*_{t,\alp}$ denote the formal adjoint of $\n_{t,\alp}$
with respect to the metrics $g^E, h^\F$. We denote by
$\Del_{t,\alp}^p$  the restriction of $\Del_{t,\alp}$
on the space of $p$-forms.
\edefe
\lem{del=del}For any  $t>0,\alp>0$
  \eq{del=del}
	\Del_{t,\alp}=(\t)^{-1}\r \Del_{1,\alp} (\r)^{-1}\t.
  \end{equation}
  In particular, the operators $\Del_{t,\alp}$ and $\Del_{1,\alp}$
  have the same spectrum.
\elem
\prf
  Clearly, the operator  $(\t)^{-1}\r$ is orthogonal, i.e.
  the adjoint of  $(\t)^{-1}\r$ is equal to
  $(\r)^{-1}\t$. The lemma follows now from  \refl{de-de'}.
\epr

\ssec{wit}{The spectrum of $\Del_{1,\alp}$}
A simple calculation \cite[Proposition 11.13]{cfks} shows that
\eq{lapl'}
	\Del_{1,\alp}=\Del+\alp A+\alp^2|df|^2,
\end{equation}
where  $\Del=\Del_{1,0}$ is the usual
Laplacian associated with the metrics $g^E, h^\F$ and $A$ is a zero order
operator.

Using the theory of globally elliptic differential operators
\cite[Ch. IV]{sh1}, one easily obtains that, for $\alp>0$ large enough,
$\Del_{1,\alp}$
has a discrete spectrum, and  that the corresponding eigenspaces in \we{}\
have finite dimension.

The following theorem computes explicitly the cohomology of the deformed
differential $\n_{1,\alp}$ on the space of $L_2$-forms.

\th{bismut}Let $m$ denote the fiber dimension of $E^-$ and let $o$
  denote the orientation bundle of $E^-$.
  If $\alp>0$ is large enough, then
  \eq{bismut}
	\dim\Ker\Del^p_{1,\alp}= \dim H^{p-m}(Z,\F_{|_Z}\otimes o)
  \end{equation}
  for any $p=0,1\nek n$.
\eth

The proof is given in the next section. In the case where $\F$ is a
trivial line bundle, \reft{bismut} was established by Bismut
\cite[Theorem 2.13]{bis}.  Note that our proof is rather different
from \cite{bis}.


\sec{L2cohom}{Proof of \reft{bismut}}

In this section we use the notation of \refs{bundle}.

\ssec{L2bis}{A cohomological interpretation of \protect$\protect\Ker
\Del_{1,\alp}$}
Assume that $\alp>0$ is sufficiently large so that the operator
$\Del_{1,\alp}$ has a discrete spectrum (cf. \refss{wit}).

Let $\wess{\b}\i\we{\b}$ be the space of smooth ($C^\infty$) square
integrable forms $\bet\in \we{\b}$ having the property
$\n_{1,\alp}\bet\in \we{\b}$.  We denote by
$H^\b_{(2)}(E,\F;\n_{1,\alp})$ the cohomology of the complex
$$
	   0\to\wess{0}@>{\n_{1,\alp}}>> \wess{1}@>{\n_{1,\alp}}>>
			\cdots @>{\n_{1,\alp}}>>\wess{n}\to 0.
$$
One should think of $H^\b_{(2)}(E,\F;\n_{1,\alp})$ as of ``deformed
$L_2$ cohomology'' of $E$ with coefficients in $\F$ (recall from
\refss{bismut} that
$\n_{1,\alp}=e^{-\alp f}\n e^{\alp f}=\n+\alp df$).
\lem{ker=coh}For any $p=0,1\nek n$, the following equality holds
  \eq{ker=coh}
	\dim \Ker\Del^p_{1,\alp}= \dim H^p_{(2)}(E,\F;\n_{1,\alp}).
  \end{equation}
\elem
\prf
Let $\widetilde{\n}_{1,\alp}, \widetilde{\n}^*_{1,\alp}$ denote the
restriction of the operators $\n_{1,\alp}, \n^*_{1,\alp}$ on
\wess{\b}. To prove the lemma it is enough to show that the following
decomposition holds
\eq{tilHo}
	\wess{\b}=\Ker\Del_{1,\alp}\oplus \IM\widetilde{\n}_{1,\alp}
	 	\oplus \IM\widetilde{\n}^*_{1,\alp}.
\end{equation}
Since the spectrum of $\Del_{1,\alp}$ is discrete, the space
$\we{\b}$ decomposes into an orthogonal  direct sum of closed subspaces
$$
	\we{\b}=\Ker\Del_{1,\alp}\oplus \IM\Del_{1,\alp}.
$$
Moreover, $\Ker\Del_{1,\alp}\i\wess{\b}$. Hence, any $\bet\in
\wess{\b}$ may be written as
\eq{bet=...}
	\bet=\kap+\Del_{1,\alp}\vphi=
	   \kap+\n_{1,\alp}\n_{1,\alp}^*\vphi+\n_{1,\alp}^*\n_{1,\alp}\vphi
\end{equation}
where $\kap\in \Ker\Del_{1,\alp}$ and $\vphi\in \we{\b}$. Then
$\Del_{1,\alp}\vphi\in \wess{\b}$.

Since $\Del_{1,\alp}$ is an elliptic operator, the form $\vphi$ is
smooth. Let $\langle\cdot,\cdot\rangle$ denote the $L_2$ scalar
product on \we{\b}.  For $\sig\in \we{\b}$, we denote by $\|\sig\|$ its
norm with respect to this scalar product. Then
$$
	\| \n^*_{1,\alp}\vphi\|^2+ \| \n_{1,\alp}\vphi\|^2=
		\langle \Del_{1,\alp}\vphi,\vphi\rangle <\infty.
$$
Hence, $\n^*_{1,\alp}\vphi,\n_{1,\alp}\vphi\in \we{\b}$. Since the
form $\vphi$ is smooth, so are the forms $\n^*_{1,\alp}\vphi,
\n_{1,\alp}\vphi$. Also, $\n_{1,\alp}\n_{1,\alp}\vphi=0$ and, by
\refe{bet=...},  $\n_{1,\alp}\n_{1,\alp}^*\vphi\in \we{\b}$. Whence
$\n^*_{1,\alp}\vphi,\n_{1,\alp}\vphi\in \wess{\b}$.
The
decomposition \refe{tilHo} follows now from \refe{bet=...}.
\epr

In view of \refl{ker=coh}, to prove the \reft{bismut} we only need to
show that
$$
	\dim H^\b_{(2)}(E,\F;\n_{1,\alp})= \dim H^{\b-\m}(Z,\F_{|_Z}\otimes o).
$$
(Recall that $\m$ denote the dimension of fibers of $E^-$).  In
fact, we will prove that the complexes
$\big(\wess{\b},\n_{1,\alp}\big)$ and
$\big(\Ome^\b(Z,\F_{|_Z}),\n\big)$ are homotopy equivalent.
\ssec{Thom}{The Thom isomorphism}Let $U^-$ be a {\em Thom form} of the
bundle $E^-$, i.e. a closed  {\em compactly supported} differential $\m$-form
on $E^-$ with values in the orientation bundle $o$ of $E^-$ whose
integrals over the fibers of $E^-$ equal 1.

We denote by $ \Ome_c^\b(E^-,\F_{|_{E^-}})$ the space of compactly
supported differential forms on $E^-$ with values in the restriction
$\F_{|_{E^-}}$ of the bundle $\F$ on $E^-$.
It is well known that the map
\eq{Thom}
	\Ome^\b(Z,\F_{|_Z}\otimes o)\to \Ome_c^{\b+\m}(E^-,\F_{|_{E^-}}),
		\qquad \bet\mapsto U^-\wedge \bet,
\end{equation}
induces the {\em Thom isomorphism}
\eq{th-is}
	H^\b(Z,\F_{|_Z}\otimes o)\to H_c^{\b+\m}(E^-,\F_{|_{E^-}})
\end{equation}
from the cohomology of $Z$ with coefficients in the bundle
$\F_{|_Z}\otimes o$
onto the compactly supported cohomology of $E^-$ with coefficients in
$\F_{|_{E^-}}$.

Let $p:E^-\to Z$ denote the natural projection. Using the
trivialization of $\F$ along the fibers of $E^-$ determined by the flat
connection $\n$ one can define the push-forward map (``integration
along the fibers'')
$$
	p_*:\Ome_c^\b(E^-,\F_{|_{E^-}})\to
		\Ome^{\b-\m}(Z,\F_{|_Z}\otimes o).
$$
The map $p_*$ induces the map
$$
	 H_c^\b(E^-,\F_{|_{E^-}})\to H^{\b-\m}(Z,\F_{|_Z}\otimes o),
$$
which is inverse to the Thom isomorphism.
\rem{p}Note that the map $p_*$ may be extended from
$\Ome_c^\b(E^-,\F_{|_{E^-}})$ to a wider class of ``rapidly decreasing
forms''. In particular, if $\bet\in \Ome^\b(E^-,\F_{|_{E^-}})$ is a
square integrable form, then the form $p_*(e^{\alp f}\bet)$ is well
defined. This remark will be used later.
\erem

\ssec{maps}{Maps between $\tilOme_{(2)}^{\b}(E,\F)$
	and $\Ome^\b(Z,\F_{|_Z}\otimes o)$}
Let $j:E^- \to E$ be the inclusion and let
$$
	j^*:\Ome^\b(E,\F)\to \Ome^\b(E^-,\F_{|_{E^-}})
$$
be the corresponding ``pull-back'' map. Set
\eq{phi}
	\phi:\wess{\b}\to \Ome^{\b-\m}(Z,\F_{|_Z}\otimes o);
		\qquad \phi:\bet\mapsto p_*j^* (e^{\alp f}\bet).
\end{equation}
By  \refr{p}, this map is well defined. Clearly, $\phi \n_{1,\alp}=
\n'\phi$, where $\n'$ is the covariant derivative on the bundle
$\F_{|_Z}\otimes o$. Hence, $\phi$ induces a map
\eq{phi-ast}
	\phi_*: H^\b_{(2)}(E,\F;\n_{1,\alp})\to
				H^{\b-\m}(Z,\F_{|_Z}\otimes o).
\end{equation}
Denote by  $\pi:E\to E^-$ the natural projection and let
$$
	\pi^*:\Ome^\b(E^-,\F_{|_{E^-}})\to \Ome^\b(E,\F)
$$
be the corresponding ``pull-back''  map defined by means of the
connection $\n$ (cf. \refss{dE'}). Set
\eq{psi}
	\psi:\Ome^{\b}(Z,\F_{|_Z}\otimes o)\to \wess{\b+\m};
		\qquad  \psi:\del\mapsto e^{-\alp f}\pi^*(U^-\wedge p^*\del).
\end{equation}
Then $\psi \n'=\n_{1,\alp}\psi$. Hence, $\psi$ induces the map
\eq{psi-ast}
	\psi_*: H^\b(Z,\F_{|_Z}\otimes o)\to H^{\b+\m}_{(2)}(E,\F,\n_{1,\alp}).
\end{equation}
\prop{homot}The maps $\phi_*$ and $\psi_*$ are inverse of each other.
\eprop

\noindent
{\it Proof.} \ Clearly, $\phi \psi=id$, and, hence, $\phi_* \psi_*=id$. We
 will show now that the map  $\psi\phi$ is homotopic to $id$.

 Consider the action of the group $\RR^*$ of nonzero real numbers on
$E=E^+\oplus E^-$ defined by
\eq{action}
	h_t:y=(y^+,y^-)\mapsto (ty^+,y^-), \qquad t\in\RR^*.
\end{equation}
As in \refss{dE'}, the flat connection on $\F$ defines naturally
the pull-back map
$$
	h_t^*:\wess{\b}\to \wess{\b}
$$
associated with $h_t$.

Assume that $\calR$ is the vector field on $E$ generating the action
\refe{action} (the {\em Euler vector field} in the direction $E^+$).
Let $\iota(\calR)$ denote the interior multiplication by $\calR$.
If
\eq{Cartan}
	\calL(\calR)=\n\iota(\calR) +\iota(\calR)\n
\end{equation}
denote the {\em Lie derivative} along $\calR$, then
\eq{Euler}
	\frac{d}{dt}h_t^*(\bet)=h_t^*(\calL(\calR)\bet)t^{-1}.
\end{equation}

Define a map $H:\wess{\b}\to \wess{\b-1}$ by
\eq{homot}
	H\bet=\int^1_0 h_t^*(\iota(\calR)\bet)t^{-1}dt,
					\qquad \bet\in\wess{\b},
\end{equation}
and let $\tilH=e^{-\alp f}He^{\alp f}$.
Note that the integral in the definition of $H$ converges because
$\calR$ vanishes at $E^-$.

\lem{homot}The following homotopy formula holds
\eq{homot'}
	\bet-\psi(\phi \bet)=(\n_{1,\alp}\tilH+\tilH \n_{1,\alp})\bet,
		\quad {\text for\  all }\quad \bet\in \wess{\b}.
\end{equation}
\elem
\prf
  If $\bet\in \wess{\b}$, let $\bet_t=h_t^*\bet$,
  and observe that $\bet_0=e^{\alp f}\psi(\phi e^{-\alp f} \bet)$ and
  $\bet_1=\bet$. Differentiating by $t$, using
  \refe{Cartan},\refe{Euler} and integrating we get
  $$
	\bet-e^{\alp f}\psi(\phi e^{-\alp f} \bet)=
			(\n H+H\n)\bet,
		\quad \text{ for  all }\quad \bet\in \wess{\b}.
  $$
  Using $\n_{1,\alp}=e^{-\alp f}\n e^{\alp f}$, we obtain the lemma.
\epr

\refl{homot}, combined with the relation $\phi\psi=id$, completes the
proof of  \refp{homot} and, hence, of \reft{bismut}.

\sec{manif}{The Bismut Laplacian on a manifold}

Let \w\ be a closed real valued differential 1-form on a compact
manifold $M$, which is non-degenerate in the sense of Bott
(cf. \refss{main}).  In this section we construct a two-parameter
family $\Del_{t,\alp}$ of Laplacians acting on the space
$\Ome^\b(M,\F)$ of differential forms on a compact manifold $M$ with
values in a flat complex vector bundle $\F$. This construction is
similar to the construction of Bismut \cite{bis}. Next we show that,
for large values of the parameter $t$, the $\dim\Ker\Del_{t,\alp}$ is
equal to the number of ``small'' eigenvalues of the Bismut Laplacian
acting on the normal bundle to the zero set of \w, where the normal
bundle is considered as a non-compact manifold.

\ssec{nuDel}{The Bismut Laplacian on the normal bundle} Throughout
this section we use the notations introduced in \refss{main}. In
particular, $C$ denotes the set of critical points of \w, i.e. the
subset of $M$ on which \w\ vanishes.  Recall that the form \w\ is
assumed to be non-degenerate in the sense of Bott, i.e.  $C$ is a
union of disjoint compact connected manifolds and the Hessian of \w\
is a non-degenerate quadratic form on the normal bundle $\nu(C)$ to
$C$ in $M$.

As we have mentioned in \refss{main}, the bundle $\nu(C)$ splits into
the Whitney sum of two subbundles
\eq{split'}
	\nu(C)=\nu^+(C)\oplus \nu^-(C)
\end{equation}
such that the Hessian of \w\ is strictly positive on $\nu^+(C)$ and
strictly negative on $\nu^-(C)$.

Denote by $p:\nu(C)\to C$ the natural projection and by
$\F_{|_{C}}$ the restriction of the bundle $\F$ on $C$. Then the
pull-back $p^*\F_{|_{C}}$ is a flat vector bundle over $\nu(C)$ which
will be denoted by $\F_\nu$.

For each connected component $Z$ of $C$ we denote by
$\nu(Z)=\nu(C)_{|_Z}$ the normal bundle to $Z$ in $M$ and by $\F_Z$ the
restriction of the bundle $\F_\nu$ on $\nu(Z)$.  Let $g^{\nu(Z)}$ be a
Riemannian metric on the non-compact manifold $\nu(Z)$ constructed as
explained in \refss{metric}, starting from a Riemannian metric on $Z$,
an Euclidean metric on $\nu(C)$ and an Euclidean connection on $\nu(C)$.
Let $h^{\F_Z}$ be a Hermitian metric on $\F_Z$
constructed as explained in \refss{metF}. We denote by $\Del^{Z}_{t,\alp}$ the
Bismut Laplacian associated with the metrics $g^{\nu(Z)}, h^{\F_Z}$
(cf. \refd{lapl}).

We denote by $g^{\nu(C)}$  the Riemannian metric on $\nu(C)$ induced by
the metrics $g^{\nu(Z)}$ and by $h^{\F_\nu}$  the Hermitian
metric on $\F_\nu$ induced by the metrics $h^{\F_Z}$.

Let
$$
	\nd=\bigoplus_{Z} \Del^{Z}_{t,\alp}:
	\Ome^\b\big(\nu(C),\F_\nu\big)\to \Ome^\b\big(\nu(C),\F_\nu\big)
$$
(the sum is taken over all connected components of $C$)
be the operator whose restriction on $\Ome^\b\big(\nu(Z),\F_Z\big)$
equals $\Del^{Z}_{t,\alp}$. Clearly,
$$
	\Ker \nd=\bigoplus_{Z} \Ker\Del^{Z}_{t,\alp}.
$$

\ssec{metrM}{Metrics on $M$ and $\F$}By the generalized Morse lemma
\cite[Ch. 6]{h} there exist  a neighborhood
$U$ of the zero section in $\nu(C)$ and an embedding
$\psi:U\to M$ such that the restriction of $\psi$ on $C$ is the
identity map and if $y=(y^+,y^-)\in U$, then
\eq{y2}
	\big(f\circ \psi\big) (y)=\frac{|y^+|^2}2-\frac{|y^-|^2}2.
\end{equation}
In the sequel, we will identify $U$ and $\psi(U)$. In particular,
we will consider $g^{\nu(C)}$ as a metric on $\psi(U)$ and
$h^{\F_\nu}$ as a metric on the restriction $\F_{|_U}$ of $\F$ on $U$.
Let $g^M$ be any Riemannian metric on $M$ whose restriction on
$U$ equals $g^{\nu(C)}$ and let $h^\F$ be any Hermitian metric on $\F$
whose restriction on $U$ equals   $h^{\F_\nu}$.

\ssec{laplM}{The Bismut Laplacian on $M$} Let $V$ be a neighborhood of
$C$ in $\nu(C)$ whose closure is contained in $U$.  Fix a function
$\phi\in C^\infty(\nu(C))$ such that $0\le \phi\le 1$, $\phi(x)=1$ \/ if
\/ $x\in V$ \/ and \/ $\phi(x)=0$ \/ if \/ $x\not\in U$.

Recall that the map $\t$ was defined in \refss{bigrad}. Using the
identification $\psi:U\to \psi(U)$ we can consider $\phi$ as a function on $M$
and $\phi\t$ as an operator on $\Ome^\b(M,\F)$.

For any $t>0,\alp>0$, we define a new differential $\n_{t,\alp}$ on
$\Ome^\b(M,\F)$ by the formula
\eq{dM}
	\n_{t,\alp}=(\phi\t+(1-\phi))^{-1}
		(\n+t\alp\w)(\phi\t+(1-\phi)).
\end{equation}
Here $\n+t\alp\w:\Ome^\b(M,\F)\to \Ome^{\b+1}(M,\F)$ denotes the map
$$
	\tet\mapsto \n\tet+t\alp\w\wedge\tet, \qquad \tet\in \Ome^\b(M,\F).
$$
\rem{not=bis}Our definition of $\n_{t,\alp}$ is slightly different from
     \cite{bis} even in the case when \w\ is cohomologious to 0 and $\F$
      is a trivial line bundle.
      However this  difference does not change the asymptotic behavior
      of the spectrum of the corresponding Laplacian (cf. \refd{DelM}).
\erem
\rem{dM=dE}
     Note that on $V$ the formulae \refe{dM} and \refe{dE} coincide.
     Also on $M\backslash U$ we have $\n_{t,\alp}=\n+t\alp \w$.
\erem
\defe{DelM}For $t>0,\alp>0$, the {\em Bismut  Laplacian} on $M$ is
   the operator
   \eq{DelM}
	\Del_{t,\alp}=\n_{t,\alp}\n^*_{t,\alp}+\n^*_{t,\alp}\n_{t,\alp}:
		\Ome^\b(M,\F)\to \Ome^\b(M,\F),
   \end{equation}
   where $\n_{t,\alp}^*$ is the formal adjoint of $\n_{t,\alp}$ with
   respect to the metrics $g^M, h^\F$.  We denote by
   $\Del_{t,\alp}^p$  the restriction of $\Del_{t,\alp}$
   on $\Ome^p(M,\F)$.
\edefe

\ssec{N}{}
We now fix a number $\alp>0$ large enough, so that the spectrum of
$\Del^C_{1,\alp}$ is discrete and the  equality \refe{bismut} holds.

Let $A$ be a self-adjoint operator with discrete spectrum.
For any $\lam>0$, we denote by $N(\lam,A)$ the number of the eigenvalues of
$A$  not exceeding $\lam$ (counting multiplicity).

The following theorem plays a central role in our proof of \reft{main}.
\th{spectr}Let $\lam_p$ $(p=0,1\nek n)$ be the smallest
   non-zero eigenvalue of \ndp. Then for any $\eps>0$ there exists
   $T>0$ such that for all $t>T$
   \eq{DelE=M}
	N(\lam_p-\eps,\Del^p_{t,\alp})=\dim\Ker\ndp.
   \end{equation}
\eth
The rest of this section is occupied with the proof of \reft{spectr}.
\ssec{above}{Estimate from above on $N(\lam_p-\eps,\Del^p_{t,\alp})$}
We will first show that
\eq{above}
	N(\lam_p-\eps,\Del^p_{t,\alp})\le\dim\Ker\ndp.
\end{equation}
To this end we will estimate the operator \H\ from below. We will
use the technique of   \cite{sh2},  adding some necessary
modifications.

Recall from \refss{metrM} that $U$ is a tubular neighborhood of
the zero section in $\nu(C)$ and that we have fixed an embedding
$\psi:U\to M$. Also in \refss{laplM} we have chosen $V\i U$.

For  $x\in U$, we denote by $|x|$ its norm with respect to the fixed Euclidean
structure on $\nu(C)$.

There exists $\kap>0$
such that the set $\{x\in U:\, |x|<2\kap\}$ is contained in $V$.
Let us fix a  $C^\infty$ function $j:[0,+\infty)\to [0,1]$ such that
$j(s)=1$ \/ for \/  $s\le \kap$, $j(s)=0$ \/ for \/  $s\ge 2\kap$
and the function $(1-j^2)^{1/2}$ is $C^\infty$. We define functions
$J,\oJ\in C^\infty(\nu(C))$ by
$$
	J(x)=j(|x|); \qquad \oJ(x)=\big(1-j(|x|)^2\big)^{\frac12}.
$$
Using the diffeomorphism $\psi:U\to \psi(U)$ we can and we will
consider $J,\oJ$ as functions on $M$.

We identify the functions $J,\oJ$ with the corresponding multiplication
operators.
For  operators  $A,B$, we denote by $[A,B]=AB-BA$ their commutator.

The following version of IMS localization formula (cf. \cite{cfks})
is due to Shubin \cite[Lemma 3.1]{sh2}.
\lem{sh1}The following operator identity holds
  \eq{sh1}
	\H=\oJ \H\oJ+J\H J+\frac12[\oJ,[\oJ,\H]]+\frac12[J,[J,\H]].
  \end{equation}
\elem
\prf
  Using the equality $J^2+\oJ^2=1$ we can write
  $$
	\H=J^2\H+\oJ^2\H=J\H J +\oJ\H\oJ+J[J,\H]+\oJ[\oJ,\H].
  $$
   Similarly,
  $$
	\H=\H J^2+\H\oJ^2=J\H J +\oJ\H\oJ-[J,\H]J-[\oJ,\H]\oJ.
  $$
  Summing these identities and dividing by 2 we come to \refe{sh1}.
\epr
We will now estimate each one of the summands in the right hand side
of \refe{sh1}.
\lem{sh3}
  There exist $c>0,\ T>0$ such that, for any $t>T$,
  \eq{sh3}
	\oJ\H\oJ\ge ct\oJ^2 I.
  \end{equation}
\elem
\prf
   Let $\langle\cdot,\cdot\rangle$ be the natural $L_2$ scalar product
   on $\Ome^\b(M,\F)$ determined by the metrics $g^M, h^\F$. For $x\in
   \Ome^\b(M,\F)$, we
   denote by $\|x\|$ its norm with respect to this scalar product.

   We can assume that $T>1$.
   Using \refe{dM}, one easily checks that for any $a\in \Ome^p(M,\F)$
   and any $t>T>1$
   $$
	\|\n_{t,\alp}a\|^2\ge
		\frac 1t\|(\n+t\alp\w) a\|^2,
	\qquad
	\|\n^*_{t,\alp}a\|^2\ge
		\frac 1t\|(\n+t\alp\w)^* a\|^2.
   $$
   Hence,
   \begin{multline}\label{E:sh3''}
	\langle \oJ\H\oJ a,a\rangle =
 	  \| \n_{t,\alp}\oJ a\|^2+\|\n_{t,\alp}^*\oJ a\|^2 \\
	   \ge \frac 1t\|(\n+t\alp\w)\oJ a\|^2+
	    \frac 1t\|(\n+t\alp\w)^*\oJ a\|^2=
	     \frac 1t\langle \oJ\Del^p_{1,\alp t}\oJ a,a\rangle.
   \end{multline}
   A simple calculation (cf. \cite[Proposition 11.13]{cfks}) shows that
   \eq{lapl''}
	\Del_{1,t\alp}=\Del+\alp t A+\alp^2t^2|\w|^2,
   \end{equation}
   where $A$ is a zero order differential operator
   and $\Del=\n\n^*+\n^*\n$ is the undeformed
   Laplacian associated with the metrics $g^M, h^\F$.

   From \refe{lapl''} and \refe{sh3''}, we get
   \eq{sh3'}
	\oJ\H\oJ\ge \frac 1t \oJ\Del^p_{1,\alp t}\oJ=
		\frac1t\oJ\Del^p\oJ+\alp \oJ A\oJ+
			t\alp^2\oJ|\w|^2\oJ.
   \end{equation}
   Since $A$ is a zero order operator, there exists $M>0$
   such that  $A>-M$. Also $\Del^p\ge 0$. Taking
   $$
	c=\frac{\alp^2}2\min_{x\in \supp\oJ}|\w|^2, \qquad T=\frac{M\alp}{c},
   $$
   and using \refe{sh3'} we get \refe{sh3}
\epr

Let $P^p_{t,\alp}:\Ome_{(2)}\big(\nu^-(C),\F_\nu\big)\to \Ker\ndel{p}{t}$
be the orthogonal projection.
This is a finite rank operator on $\Ome_{(2)}\big(\nu^-(C),\F_\nu\big)$ and its
rank equals $\dim\Ker\ndp$. Clearly,
\eq{ge}
	\ndel{p}{t}+P^p_{t,\alp}\ge \lam_p I.
\end{equation}
Using the identification $\psi:U\to \psi(U)$ we can consider
$JP^p_{t,\alp}J$  and $J\ndel{p}{t}J$ as  operators on $\Ome^p(M,\F)$.
It follows from \refr{dM=dE} that $J\H J=J\ndel{p}{t}J$.
Hence, \refe{ge} implies the following
\lem{local}
  For any $t>0$
  \eq{local}
	J\H J+JP^p_{t,\alp}J\ge \lam_pJ^2I, \qquad
		\rk JP^p_{t,\alp}J\le \dim\Ker\ndp.
  \end{equation}
\elem
For an operator $A:\Ome^p(M,\F)\to \Ome^p(M,\F)$, we denote by $\|A\|$ its norm
with
respect to $L_2$ scalar product on $\Ome^p(M,\F)$.
\lem{sh2}There exists $C>0$ such that
  \eq{sh2}
	\|[J,[J,\H]\|\le Ct^{-1}
  \end{equation}
  for any $t>0, \alp>0$.
\elem
\prf
  The left hand side of \refe{sh2} is a  zero order operator
  supported on $V$. Hence, it is enough to estimate
  $[J,[J,\ndel{p}{t}]]$.  Since
  $$
	(\r)^{-1}J(x)\r=J(xt^{-\frac12}), \qquad
		(\t)J(x)(\t)^{-1}=J(x),
  $$
  \refl{del=del} implies
  \eq{sh2'}
	[J,[J,\ndel{p}{t}]]=(\t)^{-1}\r
		[J(xt^{-\frac12}),
		  [J(xt^{-\frac12}),\ndp]]
			(\r)^{-1}\t.
  \end{equation}
  Since the operators \/ $(\t)^{-1}\r$ \/ and \/ $(\r)^{-1}\t$ \/ are mutually
  adjoint,  \refe{sh2'} implies
  $$
	\|[J,[J,\ndel{p}{t}]]\|=
		\|[J(xt^{-\frac12}),[J(xt^{-\frac12}),\ndp]]\|.
  $$
  Recall that $g^{\nu(C)}$ denotes the Riemannian metric on $\nu(C)$.
  Suppose that we have chosen local coordinates $(x^1\nek x^n)$
  near a point $x\in \nu(C)$. Set
  $$
	g_{i,j}=g^{\nu(C)}(\d/\d x^i,\d/\d x^j), \qquad 1\le i,j\le n.
  $$
  Using  \refe{lapl'}, we see that
  $$
 	\ndp=-\sum_{i,j=1}^n
		g_{i,j}	\frac{\d^2}{\d x^i \d x^j}+B,
  $$
  where  $B$ is an order 1 differential operator. Hence,
  $$
	[J(xt^{-\frac12}),[J(xt^{-\frac12}),\ndp]]=-2\sum_{i,j=1}^n g_{i,j}
		\frac{\d J(xt^{-\frac12})}{\d x^i}
			\frac{\d J(xt^{-\frac12})}{\d x^j}.
  $$
  Since the  derivatives of
  $J(xt^{-\frac12})=j(|x|t^{-\frac12})$
  can be estimated as $O(t^{-\frac12})$ the
  inequality \refe{sh2} follows immediately.
\epr
Similarly, one shows that
\eq{sh2oJ}
	 \|[\oJ,[\oJ,\H]\|\le Ct^{-1}.
\end{equation}
 From \refl{sh1}, \refl{sh3}, \refl{local}, \refl{sh2} and
\refe{sh2oJ} we get the following
\cor{above}
   For any $\eps>0$, there exists $T>0$ such that for any $t>T$
   \eq{a'}
	\H+JP^p_{t,\alp}J\ge (\lam_p-\eps)I, \qquad
		\rk JP^p_{t,\alp}J \le \dim\Ker\ndp.
   \end{equation}
\ecor
The estimate \refe{above} follows now from \refc{above} and the
following general lemma \cite[p. 270]{rs}.
\lem{general}
   Assume that $A, B$ are self-adjoint operators in a Hilbert space
   $\calH$ such that  $\rk B\le k$ and there exists $\mu>0$
   such that
   $$
	\langle (A+B)u,u\rangle \ge \mu\langle u,u\rangle
		\quad \text{for any} \quad u\in\Dom(A).
   $$
   Then $N(\mu -\eps, A)\le k$ for any $\eps>0$.
\elem

\ssec{below}{Estimate from below on $N(\lam_p-\eps, \Del^p_{t,\alp})$}
To prove \reft{spectr} it remains to  show now that
\eq{below}
	N(\lam_p-\eps,\H)\ge \dim\Ker\ndp.
\end{equation}
Let $E^p_{t,\alp}$ be the subspace of $\Ome^p(M,\F)$ spanned by the
eigenvectors of $\Del^p_{t,\alp}$ corresponding to the eigenvalues
$\lam\le \lam_p-\eps$ and let
$\Pi^p_{t,\alp}:\Ome^p(M,\F)\to E_{t,\alp}^p$ be the orthogonal
projection. Then
\eq{rkPi}
	\rk \Pi^p_{t,\alp}=N(\lam_p-\eps,\Del^p_{t,\alp}).
\end{equation}
Using the diffeomorphism $\psi:U\to \psi(U)$
we can consider $J\Pi^p_{t,\alp}J$ as an operator on
$\Ome_{(2)}\big(\nu^-(C),\F_\nu\big)$.
The proof of the following lemma does not
differ from the proof of \refc{above}.
\lem{below} For any $\del>\eps$, there exists $T>0$ such that for any
   $t>T$
   \eq{a''}
	\Del^{C,p}_{t,\alp}+J\Pi^p_{t,\alp}J\ge (\lam_p-\del)I.
   \end{equation}
\elem
The estimate \refe{below} follows now from \refe{rkPi}, \refl{below} and
\refl{general}.

\sec{proof}{Proof of \reft{main}}

Recall that the number $\alp>0$ was fixed in \refss{N}. We will use
the notation of \refs{manif}. In particular, $\lam_p \ (p=0,1,2,\dots)$
denotes the smallest non-zero eigenvalues of $\Del^{C,p}_{1,\alp}$.
By \reft{spectr} and \refl{1.2}, we can choose $t>0$ such
that
\eq{N=ker}
	N(\lam_p/2, \Del^p_{t,\alp})=\dim\Ker \Del^{C,p}_{1,\alp},
		\qquad p=0,1\nek n,
\end{equation}
and the Novikov numbers $\bet_p(\xi,\F)$ are equal to the dimension of
the cohomology of the de Rham complex $\Ome^\b(M,\F)$ with the
deformed differential
$$
	\n_t\tet=\n\tet+t\w\wedge\tet, \qquad \tet\in \Ome^\b(M,\F).
$$
The latter complex is isomorphic to the complex
\eq{com-ta}
	0\to \Ome^0(M,\F)@>{\n_{t,\alp}}>>\Ome^1(M,\F)@>{\n_{t,\alp}}>>\cdots
			@>{\n_{t,\alp}}>>\Ome^n(M,\F)\to 0,
\end{equation}
where $\n_{t,\alp}$ is defined by \refe{dM}.

Let $E^p_{t,\alp} \ (p=0,1\nek n)$ be the subspace of
$\Ome^p(M,\F)$ spanned by the
eigenvectors of $\Del^p_{t,\alp}$ corresponding to the eigenvalues
$\lam\le \lam_p/2$. From \refe{N=ker} and \reft{bismut}, we obtain
\eq{k=h}
	\dim E^p_{t,\alp}=
		\sum_{Z}\dim H^{p-\ind(Z)}(Z,\F_{|_Z}\otimes o(Z)),
\end{equation}
where the sum ranges over all connected components $Z$ of the set $C$
of critical points of \w.

Since the operator $\Del_{t,\alp}$
commutes with $\n_{t,\alp}$, the pair $(E^\b_{t,\alp},\n_{t,\alp})$
is a subcomplex of \refe{com-ta} and the inclusion induces
an isomorphism of cohomology
$$
	H^\b(E^\b_{t,\alp},\n_{t,\alp})=H^\b(\Ome^\b(M,\F),\n_{t,\alp})=
		H^\b(M,\rho_t\otimes \F).
$$
Hence,
\eq{h=b}
	\dim H^p(E^\b_{t,\alp},\n_{t,\alp})=\bet_p(\xi,\F), \qquad
			p=0,1\nek n.
\end{equation}
\reft{main} follows now from \refe{k=h}, \refe{h=b} by standard
arguments (cf. \cite{bott2}).



\end{document}